\documentclass[letter,10pt]{article}
\usepackage[utf8]{inputenc}
\usepackage{amssymb}
\usepackage{amsthm}
\usepackage{amsmath}
\usepackage{mathrsfs}
\usepackage{verbatim}
\usepackage{graphicx}
\usepackage{eucal}
\usepackage{enumerate}

\usepackage{geometry}
\usepackage{cite}
\usepackage{color}
\newtheorem{theorem}{Theorem}
\newtheorem{corollary}[theorem]{Corollary}

\newtheorem{proposition}[theorem]{Proposition}
\newtheorem{remark}[theorem]{Remark}
\newtheorem{definition}[theorem]{Definition}

\newcommand{\R}{\mathbb{R}}

\newcommand{\U}{\mathcal{U}}

\DeclareMathOperator{\argmax}{argmax}
\DeclareMathOperator{\id}{id}
\DeclareMathOperator{\X}{X}
\DeclareMathOperator{\Y}{Y}

\newcommand\Lone{{\mathrm{L}}^1}

\def\<#1,#2>{\left<#1,#2\right>}
\let\bar\overline

\def\PP{{\cal P}}

\newcommand{\TabTwo}[1]{ %
\begin{tabular}{@{}c@{\hspace{1mm}}c@{}}
 #1
\end{tabular}
}

\title{A note on Cournot-Nash equilibria and Optimal Transport between unequal dimensions\footnote{B.P. is pleased to acknowledge the support of National Sciences and Engineering Research Council of Canada Discovery Grant number  04658-2018.   L.N. was supported by a public grant as part of the Investissement d'avenir project, reference ANR-11-LABX-0056-LMH, LabEx LMH and  from the CNRS PEPS JCJC (2022). }}

\author{
Luca Nenna\thanks{Université Paris-Saclay, CNRS, Laboratoire de mathématiques d’Orsay, 91405, Orsay, France.  \texttt{luca.nenna@universite-paris-saclay.fr}}  
\and 
Brendan Pass\thanks{Department of Mathematical and Statistical Sciences, 632 CAB, University of Alberta, Edmonton, Alberta, Canada, T6G 2G1  \texttt{pass@ualberta.ca}} 
}

\begin{document}

\maketitle
 
 \begin{abstract}
 This note is devoted to study a class of games with a continuum of players for which Cournot-Nash equilibria can be obtained by minimisation of some cost related to Optimal Transport. In particular we focus on the case of an Optimal Transport term between unequal dimension. We also present several numerical simulations.  
 \end{abstract}

 \tableofcontents

 \section{Introduction}
 Since  Aumann's seminal works \cite{Aumann,Aumann2}, equilibria in games with a continuum of players have received a lot of attention from the Economics and Game Theory communities. Following Mas-Colell's approach \cite{MasColell}, we consider a type space $\X\subset\R^m$ endowed with a probability measure $\mu$. Each player of type $x$ has to choose a strategy $y\in\Y$ (where $\Y\subset\R^n$) in order to minimise a cost $\Phi:\X\times\Y\times\PP(\Y)\rightarrow \R$ which depends both on his type $x$ and strategy $y$ as well as the distribution of strategies $\nu\in\PP(\Y)$ resulting from the other  players' behaviour. Since the cost depends on other players' strategies only through the distribution $\nu$,  it is not important who plays a specific strategy, but how many players chose it (i.e. the game is anonymous).  Thus, a Cournot-Nash equilibrium is defined as a probability measure $\gamma\in\PP(\X\times\Y)$, whose first marginal is $\mu$, such that 
 
 \begin{equation*}
 \gamma(\{ (x,y)\in\X\times\Y\;:\;\Phi(x,y,\nu)=\min_{z\in\Y}\Phi(x,z,\nu)\})=1
 \end{equation*}     
 
where the strategy distribution $\nu$ is the second marginal of $\gamma$.
Notice that if we consider an homogenous population of players (which mean that we do not need a distribution $\mu$ anymore) then we retrieve exactly the definition of Nash equilibrium.
\\
Let us now focus on the additively separate case where the total cost $\Phi$ can be written as $\Phi(x,y,\nu)=c(x,y)+F(\nu,y)$. It has been recently showed in \cite{abgcmor} (see also \cite{abgcptrl}) that Cournot-Nash equilibria can be obtained by the minimization of a certain functional on the set of measures on the space of strategies, which actually means that the games considered in \cite{abgcmor} belongs to a class of models which have the structure of a  potential game (i.e. variational problems). This functional typically involves two terms: an optimal transport cost and a more standard integral functional which may capture both congestion and attractive effects. 
This variational approach  is somehow more constructive and informative (but less general since it requires a separable cost) than the one relying on fixed-point arguments as in \cite{MasColell} but the optimal transport  cost term is delicate to handle. It is indeed costly in general to solve numerically an optimal transport problem and compute an element of the subdifferential of the optimal cost as a function of its marginals. However it has been recently introduced (see for instance \cite{Galichon-Entropic,CuturiSinkhorn,benamouetalentropic}) a very efficient numerical method relying on an entropic regularization of optimal transport which turned to be useful also to approximate Cournot-Nash equilibria (see \cite{Blanchet2017}), in the case of 
potential games. 
\\
In \cite{abgcmor} the authors treat the case where the measures $\mu$ and $\nu$ are probabilities on $\X\subset\R^m$ and $\Y\subset\R^n$, respectively, where $m=n$. In this article we want to address the case in which $m>n$ so that the optimal transport term becomes an \emph{unequal dimensional Optimal Transport term} \cite{PassM2one,MultiMatch} .
In a recent paper \cite{NennaPass1} we have showed that for a variety of different forms of $\mathcal F[\nu](y)$ a nestedness condition, which is known to yield much improved tractability of the optimal transport problem, holds for any minimizer. Depending on the exact form of the functional, we exploit this to find local differential equations characterizing solutions, prove convergence of an iterative scheme to characterise  the solution, and prove regularity results.
The aim of this paper is to illustrate some numerical methods to compute the solution in the case of Cournot-Nash equilibria by using the characterization  results established in \cite{NennaPass1}.
\paragraph{Motivation} We want to briefly remark through an example that unequal dimensional Optimal Transport is actually quite natural in the Economics/Game theory framework as the one we are dealing with.\\
Let think that we have a population of physicians who are characterized by a type $x$ belonging to $X\subset\R^m$ where the dimension $m$ is the number of characteristics: age, gender, university where they get their diploma, their hometown, etc. 
The probability measure $\mu$ on $X$ gives then the distribution of types in the physician population.
The main idea of Optimal Transport is to match the physician with a city, for example, where they would like to open their private practice. The cities are characterized by a  type $y$ which belongs to $Y\subset \R^n$ where $n$ is now the number of characteristics and $Y$ is endowed with a probability measure $\nu$ giving the distribution of types of cities. 
 Notice that the matching between physicians and cities must minimize a given transportation cost $c(x,y)$ which can indeed model the cost for physician $x$ to commute to the city $y$ where he/she has his/her practice.
 One can also think, in more economical terms, that the cost is of the form $c(x,y)=-s(x,y)$ where $s$ is a surplus and in this case the optimal matching will maximize the given surplus.
 It is actually quite important to highlight now that the number of characteristics $n$ of the cities can be much smaller than the one for the physician; indeed it is natural to take $n=1$ saying that just a characteristic (i.e. the population) of the city is taken into account by the physicians.   
    
 \begin{remark}[Notation]
 With a slight abuse of notation we will identify the measure with its density throughout all the paper. 
\end{remark}

 \section{Optimal Transport between unequal dimensions}
 Given two probability measures $\mu$ on $X\subset\R^m$ and $\nu$ on $Y\subset\R^n$, the Monge-Kantorovich problem consists in transporting $\mu$ onto $\nu$ so as to optimise a given cost function $c(x,y)$.
 Indeed we look for a probability measure $\gamma$ on the product space $X\times Y$ whose marginals are $\mu$ and $\nu$ such that it solves the following maximisation problem 
 
 \begin{equation}
 \label{mk}
 \mathcal T_c(\mu,\nu):=\inf_{\gamma\in\Pi(\mu,\nu)}\displaystyle\int_{\X\times\Y}c(x,y)d\gamma(x,y),
 \end{equation}

 where $\Pi(\mu,\nu):=\{ \eta\in\PP(\X\times\Y)\;:\; \pi_x(\gamma)=\mu,\;\pi_{y}(\gamma)=\nu  \}$ and $\pi_x: X\times Y\rightarrow X$.
%
%

We refer the reader to \cite{Villani-OptimalTransport-09,santambook}  for results about the characterisation of optimal maps for the case in which $m=n$.
Before discussing the case $m>n$, it is important  to highlight that the Monge-Kantorovich problem  admits a dual formulation which is useful in order to understand the solution to \eqref{mk} 

\begin{equation}
\mathcal T^{dual}_c(\mu,\nu):=\sup_{(u,v)\in\U}  \displaystyle \int_{\X} u(x)d\mu(x)+\int_{\Y} v(y)d\nu(y),
\end{equation}

where $\U:=\{ (u,v)\in\Lone(\mu)\oplus\Lone(\nu)\;:\; u(x)+v(y)\geq c(x,y)\;\text{on}\;X\times Y\}$. The optimal $u,v$ are called Kantorovich potentials. The interesting fact is that $\mathcal T_c=\mathcal T^{dual}_c$.

Under mild conditions (for instance, if $c$ is differentiable, as we are assuming here, $X$ connected and $\bar \mu(x) >0$ for all $x \in X$ \cite{santambook}), there is a unique solution $(u,v)$ to the dual problem, up to the addition $(u,v)\mapsto(u+C,v-C) $ of a constant, known as the Kantorovich potentials and these potentials are $c$-concave; that is, they satisfy
$$
u(x) =v^c(x):=\min_{y \in Y}[c(x,y) -v(y)], \text{ }v(y) =u^c(y):=\min_{x \in X}[c(x,y) -u(x)],
$$
%

In particular if  both $\nu$ and $\mu$ are absolutely continuous, then the potential $v$ satisfies the Monge-Ampere type equation almost everywhere \cite{mccann2018optimal}\footnote{Note that, here and below, our notation differs somewhat from \cite{PassM2one} and \cite{mccann2018optimal}, since we have adopted the convention of minimizing, rather than maximizing, in \eqref{mk}.}:

 \begin{equation}
\label{eqn: unequal MA}
\nu(y)=\int_{\partial^cv(y)}\dfrac{\det(D_{yy}^2c(x,y)-D^2v(y))}{\sqrt{|\det(D^2_{yx}c D^2_{xy}c)(x,y)|}}\mu(x)d\mathcal H^{m-n}(x),
\end{equation}
where 
$\partial^cv(y):=\{x: u(x)+v(y) =c(x,y)\}$.  In general,  this is a \emph{non-local} differential equation for $v(y)$, since the domain of integration $\partial^cv(y)$ is defined using the values of $v$ and $u=v^c$ throughout $\Y$; however, when the model satisfies the generalized nestedness condition   (namely $\partial^cv(y) =X_=(y,Dv(y))$, where $X_=(y,p) =\{x \in X: D_yc(x,y) =p\}$, see \cite{NennaPass1}[Definition 1] for more details), it reduces to the local equation \cite{mccann2018optimal}:

 \begin{equation}
\label{eqn: local unequal MA}
\nu(y)=\int_{X_=(y,Dv(y))}\dfrac{\det(D_{yy}^2c(x,y)-D^2v(y))}{\sqrt{|\det(D^2_{yx}c D^2_{xy}c)(x,y)|}}\mu(x)d\mathcal H^{m-n}(x).
\end{equation}
\subsection{Multi-to one-dimensional optimal transport}
We consider now the optimal transport problem in the case in which $m>n=1$ (for more details we refer the reader to \cite{PassM2one}).
Let us define the super-level set  $ X_\geq (y,k)$ as follows

\begin{equation*}
X_\geq (y,k):=\{x\in\X\;:\; \partial_yc(x,y)\geq k \},
\end{equation*}

 and its strict variant $X_> (y,k):=X_\geq (y,k)\setminus X_=(y,k)$.
 In order to build an optimal transport map $T$ we take the unique level set splitting the mass proportionately with $y$; that is $k(y)$ such that
 
 \begin{equation} 
 \label{keq}
 \mu(X_\geq(y,k(y)))=\nu((-\infty,y]),    
 \end{equation} 
 
 then we set $y=T(x)$ for all $x$ which belongs to $X_=(y,k(y))$.
 Notice that if $x\in X_=(y_0,k(y_0)) \cap X_=(y_1,k(y_1))$ then the map $T$ is not well-defined. In order to avoid such a degenerate case we ask that the model $(c,\mu,\nu)$ satisfies the following property
 
\begin{definition}[Nestedness $n=1$]
\label{nested}
The model $(c,\mu,\nu)$ is nested if
\[\forall y_0,y_1\; \text{with} \; y_1>y_0,\; \nu([y_0,y_1])>0\Longrightarrow X_\geq (y_0,k(y_0)) \subset X_>(y_1,k(y_1)) .\] 
\end{definition}

 Thus, if the model $(c,\mu,\nu)$ is nested then \cite{PassM2one}[Theorem 4] assures that $\gamma_T=(\id,T)_\sharp \mu$, where the map $T$ is built as above, is the unique maximiser of \eqref{mk} on $\Pi(\mu,\nu)$. Moreover, the optimal potential $v(y)$ is given by $v(y)=\int^y k(t) dt$.
 
 We remark that  nestedness depends on all the data $(c,\mu,\nu)$ of the problem, in the following  we give a condition which can ensure nestedness when only $c$ and $\mu$ are fixed.

 Consider now the equality $\nu([-\infty,y])=\mu(X_\leq(y,k(y)))$, then differentiating it gives the following equation
 \begin{equation}
 \label{MAmulti2one}
 \nu(y)=\int_{X_=(y,k(y))}\dfrac{D_{yy}^2c(x,y)-k'(y)}{|D^2_{xy}c(x,y)|}\mu(x)d\mathcal H^{m-1}(x),
 \end{equation}
which can be viewed as the multi-to-one counterpart of the Monge-Ampère equation in classical Optimal Transport theory (note also that it is exactly the $n=1$ case of \eqref{eqn: local unequal MA}).\\
  We now briefly recall some results of \cite{NennaPass1} which guarantee the nestedness of the model by checking that a lower bound on $\nu$ satisfies a certain condition.
  
 Fix $y_0 <y_1$ (where $y_0,y_1\in\Y$), $k_0 \in D_yc(\X,y_0)$ and set $k_{max}(y_0,y_1,k_0)=\sup\{k: X_\geq(y_0,k_0) \subseteq X_\geq (y_1,k)\}$.  We then define the \emph{minimal mass difference}, $D^{min}_\mu$, as follows: 
 
 $$
D^{min}_\mu(y_0,y_1,k_0) =\mu(X_\geq(y_1,k_{max}(y_0,y_1,k_0))\setminus X_\geq(y_0,k_0)).
 $$
 The minimal mass difference represents the smallest amount of mass that can lie between $y_0$ and $y_1$, and still have the corresponding level curves $X_=(y_0,k_0)$ and $X_=(y_1,k_1)$ not intersect.
 \begin{theorem}[Thm. 3 \cite{NennaPass1}]
 \label{lemma: nestedness condition}
 	Assume that $\mu$ and $\nu$ are absolutely continuous with respect to Lebesgue measure.     If $D^{min}_\mu(y_0,y_1,k(y_0)) < \nu([y_0,y_1])$ for all $y_0<y_1$ where $k(y)$ is defined by \eqref{keq}, then $(c,\mu,\nu)$ is nested.  Conversely, if $(c,\mu,\nu)$ is nested, we must have $D^{min}_\mu(y_0,y_1,k(y_0)) \leq \nu([y_0,y_1])$ for all $y_1 > y_0$.
 \end{theorem}
 then the following condition on the density of $\nu$ follows
  \begin{corollary}[Cor. 4 \cite{NennaPass1}]
  \label{cor: sufficient conditions for nestedness}
Assume that $\mu$ and $\nu$ are absolutely continuous with respect to Lebesgue measure.  If for each $y_0 \in Y$, we have
 	$$
 	\sup_{y_1 \in Y,y_0\leq y\leq y_1}\left[\frac{D^{min}_\mu(y_0,y_1,k(y_0))}{y_1-y_0}-\nu(y)\right] < 0,
 	$$
then $(c,\mu,\nu)$ is nested.	
 \end{corollary}

 \section{A variational approach to Cournot-Nash equilibria}
We briefly recall here the variational approach to Cournot-Nash via Optimal Transport proposed  in\cite{blanchet2014remarks,abgcmor,abgcptrl}.
Agents are characterised by a type $x$ belonging to a compact metric space $X$ which is endowed with a given probability measure $\mu\in\PP(X)$ which gives the distribution of types in the agent population.
For sake of simplicity we will consider probability measures absolutely continuous with respect to the Lebesgue measure and, with a slight abuse of notation we will identify the measure with its density.
Each agent must choose a strategy $y$ from a space strategy space $Y$ (a compact metric space again) by minimizing a given cost $\Phi:X\times Y\times\PP(Y)\rightarrow \R$.
Notice that the cost $\Phi(x,y,\nu)$ of one agent does not depend only on the type of the agent and the chosen strategy but also on the other agents' choice through the probability distribution $\nu$ resulting from the whole population strategy choice.
Then, an equilibrium can be described by a joint probability distribution $\gamma$ on $X\times Y$ which gives the joint distribution of types and strategies and is consistent with the cost-minimising behaviour of agents.
 \begin{definition}[Cournot-Nash equilibrium]
 Given a cost $\Phi: X\times Y\times\PP(Y)\rightarrow \R$, a probability measure $\mu\in\PP(X)$, a Cournot-Nash equilibrium is a probability measure on $\PP(X\times Y)$ such that
 \begin{equation}
 \label{CN}
  \gamma(\{ (x,y)\in X\times Y\;:\;\Phi(x,y,\nu)=\max_{z\in Y}\Phi(x,z,\nu)\})=1,
 \end{equation}
 and $\pi_x(\gamma)=\mu$, $\pi_y(\gamma)=\nu$.
\end{definition} 
Here we consider Cournot-Nash equilibria in the separable case that is   
 \[\Phi(x,y,\nu)=c(x,y)+F(\nu,y),\]
where $c$ is a transport cost and the term $F(\nu,y)$ captures all the strategic interactions.
In particular $F$ can be of the form
\[ F(\nu,y)=f(\nu)+\int_Y\phi(y,z)d\nu(z), \]
where
\begin{itemize}
\item $f(\nu)$ captures the congestion effects: frequently played strategies are costly;
\item $\int_Y\phi(y,z)d\nu(z)$ captures the positive interactions.
\end{itemize}
\begin{remark}
In order to understand better the role played by the term $F(\nu,y)$, let us consider again the example of the physicians and the cities that we explained in the introduction.
In this case the distribution $\nu$ of the types of the cities, namely the strategies the physician have to choose, is unknown.
The main idea of the Cournot-Nash equilibria is to find the matching $\gamma$ such that the agents minimize a transport cost (again this can represent the cost of commuting) and an extra term which takes into account 
\begin{itemize}
\item a repulsive effect (the congestion): it is self-defeating opening a practice in a city chosen already by many other physicians;
\item an attractive effect  (the positive interactions term): the physicians want to share their own experiences with their peers.
\end{itemize}
\end{remark}
In \cite{abgcmor} the authors observed that $F$ is the first variation of the energy
\[ \mathcal F(\nu)=\int_Yf(\nu(y))dy+\int_{Y\times Y}\phi(y,z)d\nu(y)d\nu(z),  \]
and consequently   proved that the condition defining the equilibria is in fact the Euler-Lagrange equation for the following variational problem
\begin{equation}
\label{pbVariation}
\inf\{ \mathcal J(\mu,\nu)\;|\;\nu\in\PP(Y) \},
\end{equation}
where 
\[\mathcal J(\mu,\nu)\:=\mathcal T_c(\mu,\nu)+\mathcal F(\nu).\]
Indeed if $\nu$ solves \eqref{pbVariation} and $\gamma$ is the optimal transport plan between $\mu$ and $\nu$ for the cost $c(x,y)$ then $\gamma$ is a Cournot-Nash equilibrium for the cost $\Phi(x,y,\nu)$ defined as above.
Notice that once the cost function satisfies the twisted condition (also known in economics as the Spence-Mirrlees condition)
then $\gamma$ is determinist, that is $\gamma$ is a pure equilibrium (see, for example, \cite{santambook}).

\subsection{Nestedness of Cournot-Nash equilibrium}
We now briefly recall the main results in  \cite{NennaPass1} assuring that nestedness condition holds.
Moreover, in the next section, we will see how nestedness can lead to the development of numerical methods to solve \eqref{pbVariation}.
 
\paragraph{The congestion case} We firstly consider the case $n=1$, meaning that we are in  the multi-to-one Optimal Transport problem, and  the functional $\mathcal F(\nu)$ takes into account only the congestion effects.
Indeed we consider the following functional
\[ \mathcal F(\nu)=\int_Yf(\nu)dy\]
with $f:[0,\infty)\to\R$ continuously differentiable on $[0,\infty)$, strictly convex with superlinear growth at infity and satisfying
\[ \lim_{x\to0^+} f'(x)=-\infty.\]
An example is the entropy $f(\nu)=\nu\log(\nu)$.
In particular we have this result establishing lower and upper bounds on the density $\nu$ and the nestedness of the model.
\begin{theorem}[Thm. 11 \cite{NennaPass1}]
Assume that $\nu\in\PP(Y)$, with $Y=(0,\bar y)$ is a minimizer of \eqref{pbVariation} then $\nu$ is absolutely continuous w.r.t. the Lebesgue measure and  there exist two constants $M_1$ and $M_2$ depending on the cost function and $Y$ such that
\[(f')^{-1}(M_1)\leq \nu(y) \leq (f')^{-1}(M_1)\].
Moreover, $(c,\mu,\nu)$ is nested provided
\[ \sup_{y_1\in Y, y_0\leq y\leq y_1}\dfrac{D_\mu^{min}(y_0,y_1,k(y_0))}{y_1-y_0}-(f')^{-1}(M_1)<0 \]
for all $y_0\in Y$.
\end{theorem}
Notice that once we have  established the model is nested then by setting $k(y)=v'(y)$ we get
\begin{equation}
 \label{MAmulti2oneCong}
 \nu(y)=\int_{X_=(y,k(y))}\dfrac{D_{yy}^2c(x,y)-k'(y)}{|D^2_{xy}c(x,y)|}\mu(x)d\mathcal H^{m-1}(x):=G(y,k(y),k'(y))
 \end{equation}
 and differentiating the first order condition of \eqref{pbVariation}
 we get a second order differential equation for $k$
 \[ k(y)+f''(G(y,k(y),k'(y)))\frac{d}{dy}G(y,k(y),k'(y))=0.  \]
 In section \ref{Numerics:CC} we detail how helpful nestedness is in order to design a suitable numerical method.
\paragraph{The interaction case}
We consider now the case $m>n\geq 1$ and a functional $\mathcal F$ which captures only the interaction effects, that is
\[ F(y,\nu)=V(y)+\int_Y\phi(y,z)d\nu(z) , \]
where $\phi$ is symmetric that is $\phi(y,z)=\phi(z,y)$.
In this case it is actually more difficult to find lower and upper bounds as in the cogestion, but still  the following nestedness  result holds
\begin{theorem}[Thm. 15 \cite{NennaPass1}]
\label{thm: nestedInteraction}
Assume that $y\mapsto c(x,y)+V(y)+\phi(y,z)$ is uniform convex throughout $X\times Y \times Y$ and 
\[ \langle\; (D_yc(x,y)+DV(y)+D_y\phi(y,z))\;|\; n_Y(y)\;\rangle\geq 0\quad \forall x\in X, z\in Y, y\in\partial Y,  \]
where $n_Y(y)$ is the outward normal to $Y$.
Then the model $(c,\mu,\nu)$ is nested.
\end{theorem}
Moreover in this special case one can characterize the minimizer by using the best-reply scheme that we will detail in section \ref{Numerics:BRS}.

 \section{Numerics}
We now describe three different numerical approaches to solve \eqref{pbVariation}.
In particular we are now considering discretized measures that is $\mu=\sum_{i=1}^M\mu_i\delta_{x_i}$ ($x_i\in\R^m$) and $\nu=\sum_{i=1}^N\nu_i\delta_{y_i}$ ($y_i\in\R^n$), where $M$ and $N$ is the number of gridpoints chosen to discretize $X$ and $Y$ respectively. To reduce the amount of notation here, we use the same notations for the continuous problem as for the discretized one where integrals are replaced by finite sums and $c$ and $\gamma$ are now $M\times N$ matrices.
\subsection{Congestion case}
\label{Numerics:CC}
\begin{figure}[h!]
\TabTwo{
\includegraphics[width=.55\linewidth]{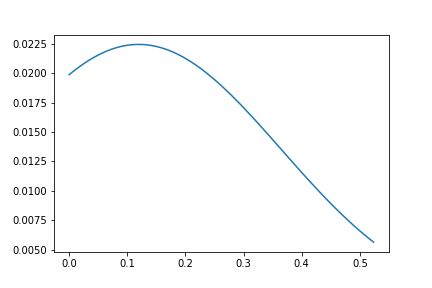}&
\includegraphics[width=.55\linewidth]{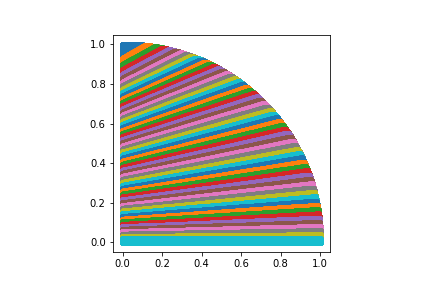}\\
}
\caption{(Left) Final density $\nu^\star$. (Right) Intersection between the level-set $X_\geq(y,k^\star(y))$ for the final solution and the support of the fix measure $\mu$}
\label{fig:quadUniform}
\end{figure}

In this section we deal with a numerical method exploiting the nestedness result we provided for the congestion case.
Consider the \[ \mathcal F(\nu)=\int_Yf(\nu)dy\]
with $f$ satisfying the hypothesis we gave in the previous section and   $Y=(0,\bar{y})$ such that the model is nested.
Then, the main idea of the iterative algorithm is actually to  combine the optimality condition with the fact that in the nestedness case we can explicitly build the optimal transport map and the Kantorovich potential.
We firstly recall that in this case the optimality condition of \eqref{pbVariation} reads
\begin{equation}
\label{OptConCong}
v(y)+f'(\nu)=C.
\end{equation}
Since the model is nested we have that $v(y)=\int_0^y k(s)ds$ which leads to 
\begin{equation} 
\label{eq:optinu}
\nu(y)=(f')^{-1}\bigg(C-\int_0^y k(s)ds\bigg) .
\end{equation}
The idea is now to use \eqref{eq:optinu} to obtain the following iterative algorithm: fix an initial density $\nu^{(0)}$
\begin{equation}
\label{eq:iterations}
\begin{cases}
&k^{(n)}(y)\;\text{s.t}\;  \mu(X_\geq(y,k^n(y)))=\nu^{(n-1)}((0,y]),\\
&\nu^{(n)}(y)=(f')^{-1}(C^{(n)}-\int_0^y k^{(n)}(s)ds), 
\end{cases}
\end{equation}
where $C^{(n)}$ is computed such that $\nu^{(n)}$ is a probability density.
\begin{remark}
This algorithm is actually an adaptation to the unequal dimensional case of the one proposed in \cite{blanchet2014remarks}. 
For both the equal and  unequal case we proved in \cite{NennaPass3} a convergence result by choosing a suitable metric.
\end{remark}
\begin{remark}
Take $f(\nu)=\nu\log(\nu)$ then the update of $\nu^{(n)}$ reads
\[ \nu^{(n)}(y)=\dfrac{\exp(-\int_0^y k^{(n)}(s)ds)}{\int_Y\exp(-\int_0^y k^{(n)}(s)ds)dy} \]
\end{remark}
\paragraph{Numerical results}
In all the simulations we present in this paragraph we have taken the uniform density on the arc $(0,\frac{\pi}{6})$ as initial density $\nu^{(0)}$, $\mu$ is always the uniform on quarter disk $X=\{x_1,x_2>0\;:\; x^2_1+x^2_2<1\}$ and $f(\nu)$ is the entropy as in the remark above. By Corollary 12 in \cite{NennaPass1}, when the cost  is given by $c(x,y)=|x_1-\cos(y)|^2+|x_2-\sin(y)|^2$ we know that the model is nested provided $\bar y\leq0.61$ (which is exactly the case we consider for the numerical tests).
Moreover, we have also considered an additional potential term, which will favour a concentration in a certain area of $(0,\frac{\pi}{6})$, so that the functional $\mathcal F$  has the form
\[ \mathcal F(\nu)=\int_Yf(\nu)dy+\int V(y)d\nu(y),\]
where $V(y)=10|y-0.1|^2$.
The second equation in \eqref{eq:iterations} then  becomes
\[ \nu^{(n)}(y)=(f')^{-1}\bigg(C^{(n)}-V(y)-\int_0^y k^{(n)}(s)ds\bigg).\]
In Figure \ref{fig:quadUniform} we show the final density we have obtained (on the left) and the intersection between the level-set $X_\geq(y,k^\star(y))$ for the final solution and the support of the fix measure $\mu$. Notice that, as expected, the level set $X_\geq(y,k^\star(y))$ do not cross each other meaning that the model is nested.

\begin{figure}[h!]
\TabTwo{
\includegraphics[width=.55\linewidth]{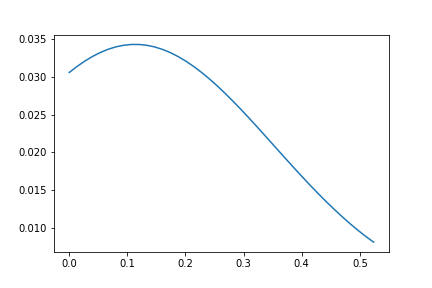}&
\includegraphics[width=.55\linewidth]{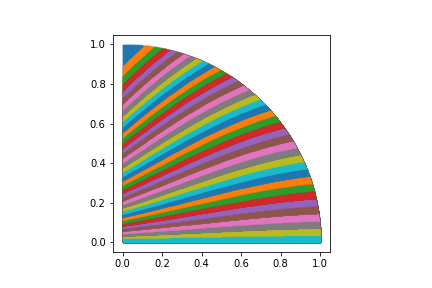}\\
\includegraphics[width=.55\linewidth]{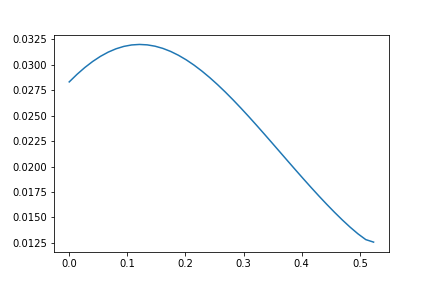}&
\includegraphics[width=.55\linewidth]{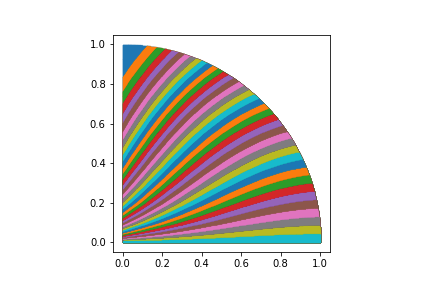}\\
}
\caption{First row: $p=4$.(Left) Final density $\nu^\star$. (Right) Intersection between the level-set $X_\geq(y,k^\star(y))$ for the final solution and the support of the fix measure $\mu$. Second row: $p=8$.(Left) Final density $\nu^\star$. (Right) Intersection between the level-set $X_\geq(y,k^\star(y))$ for the final solution and the support of the fixed measure $\mu$.}
\label{fig:pUniform}
\end{figure}

In Figure \ref{fig:pUniform}, we present some simulations by taking the same data set as above, but the cost function is now given by
\[ c(x,y)=\bigg( 1+ |x_1-\cos(y)|^2+|x_2-\sin(y)|^2 \bigg)^{\frac{p}{2}} \] 
with $p>2$. Notice that now the level set are not straight lines anymore. 

\subsection{Best reply scheme}
 \label{Numerics:BRS}
\begin{figure}[h!]
\TabTwo{
\includegraphics[width=.55\linewidth]{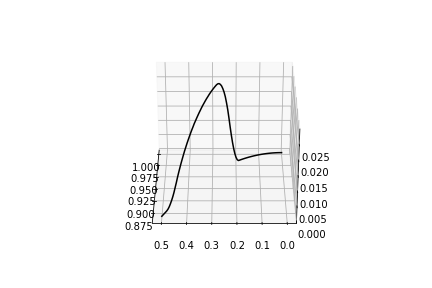}&
\includegraphics[width=.55\linewidth]{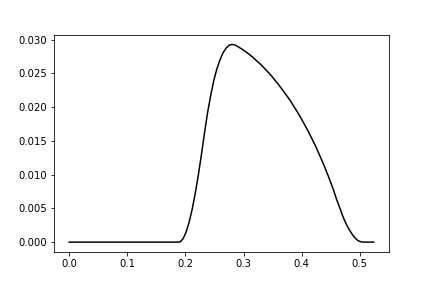}\\
}
\caption{(Left) Final density $\nu^\star$ concentrated on the arc $(0,\frac{\pi}{6})$. (Right) Final density $\nu^\star$.}
\label{fig:quadUniformBR}
\end{figure}

 Assume now that
 \[ F(y,\nu)=V(y)+\int_Y\phi(y,z)d\nu(z) \]
 and consider the first optimality condition of \eqref{pbVariation} which reads
 \begin{equation}
 \label{OptCon}
 c(x,y)+F(y,\nu)=C,
 \end{equation}
 where $C$ is a constant assuring that the minimizer $\nu$ is a probability distribution.
 Notice that by differentiating  \eqref{OptCon} with respect to $y$ we get 
\begin{equation}
\label{OPtCond}
 D_yc(x,y)+D F(y,\nu)=0 ,
 \end{equation}
 we denote by $B_\nu:\X\to\Y$ the map such that  
 \begin{equation}\label{eqn: iteration map} 
 D_yc(x,B_\nu(x))+D F(\nu,B_\nu(x))=0,
 \end{equation}
which is well defined under the hypothesis of theorem \ref{thm: nestedInteraction}.
Then, the best-reply scheme, firstly introduced in \cite{blanchet2014remarks}, consists in iterating the application defined as  
\begin{equation}
\label{It}
\mathcal B(\nu):=(B_\nu)_\sharp\mu. 
\end{equation}
 We refer the reader to \cite{NennaPass1}[Theorem 19] for a convergence result.
 Notice that \cite{NennaPass1}[Theorem 19]  assures also that $\mathcal B(\nu)$ is absolutely continuous with respect to the Lebesgue measure, meaning that the density, and so \eqref{It}, can be computed by using the co-area formula.
 We obtain then the following algorithm: fix $\nu^{(0)}$ then
 \[ \nu^{(n+1)}= \mathcal B(\nu^{(n)}):=(B_{\nu^{(n)}})_\sharp\mu.\]
 \paragraph{Numerical results}
 As in the previous section the initial density $\nu^{(0)}$ is the uniform on the arc $(0,\frac{\pi}{6})$ and $\mu$ is the uniform on $X=\{x_1,x_2>0\;:\; x^2_1+x^2_2<1\}$. Moreover, we take here a potential $V(y)=|y-\frac{\pi}{12}|^2$ and a quadratic interaction $\phi(y,z)=|y-z|^2$.
In Figure \ref{fig:quadUniformBR} we plot the final density obtained in the case in which the cost function is $c(x,y)=|x_1-\cos(y)|^2+|x_2-\sin(y)|^2$. Notice that since there is not a congestion term like the entropy, the final density has a support much smaller than the interval  $(0,\frac{\pi}{6})$.
 
\begin{figure}[h!]
\TabTwo{
\includegraphics[width=.55\linewidth]{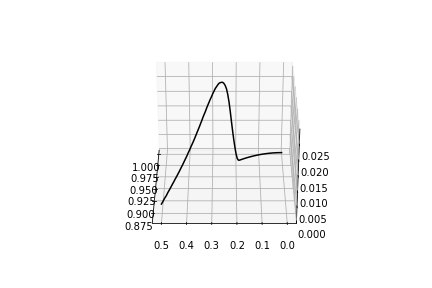}&
\includegraphics[width=.55\linewidth]{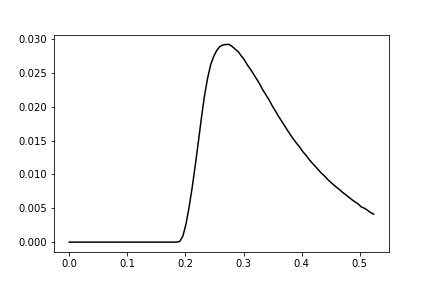}\\
\includegraphics[width=.55\linewidth]{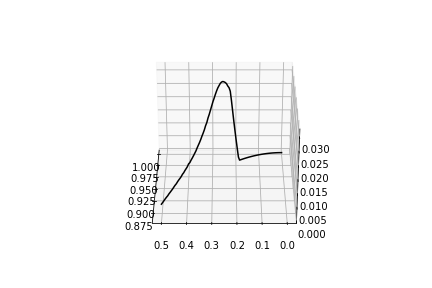}&
\includegraphics[width=.55\linewidth]{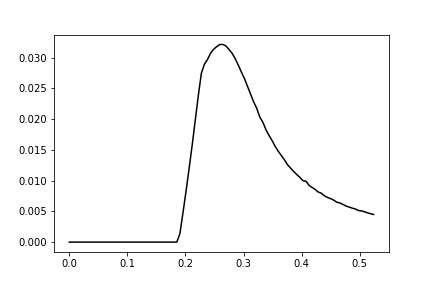}\\
}
\caption{(Left) Final density $\nu^\star$ concentrated on the arc $(0,\frac{\pi}{6})$. (Right) Final density $\nu^\star$.}
\label{fig:quadUniformBRp84}
\end{figure}

 In Figure \ref{fig:quadUniformBRp84} we repeat the simulation with the same data set but taking the cost function $c(x,y)=\bigg( 1+ |x_1-\cos(y)|^2+|x_2-\sin(y)|^2 \bigg)^{\frac{p}{2}}  $ with $p=4,8$. Notice that the cost function plays the role of an attractive term with respect to the fix measure $\mu$ which makes the final density have a more spread support. 
 
 \subsection{The Entropic Regularization}
  \label{Numerics:Ent}
 The last numerical approach we consider is the entropic reguarization of Optimal Transport (we refer the reader to \cite{benamouetalentropic,abgcmor,Chizat,CuturiSinkhorn,peyregradientflows} for more details) which can be easily used to solve \eqref{pbVariation} no matter the choice of the transportation cost $c(x,y)$ and the functional $\mathcal F(\nu)$.
 This actually implies that the entropic regularization is clearly more flexible than the methods we have previously introduced but in some cases (i.e. the congestion term is the entropy) the major drawback of this approach is an extra diffusion.
The entropic regularization of Optimal Transport can be stated as follows
\begin{equation}
\label{eq: entropic OT}
\mathcal T_{c,\epsilon}(\mu,\nu):=\inf\left\{ \int_{X\times Y}cd\gamma+\epsilon\int_{X\times Y} (\log(\gamma)-1)d\gamma\;|\;\gamma\in\Pi(\mu,\nu)  \right\},
\end{equation}
  where we assume that  $0(\log(0)-1)=0$.
  The problem now is strictly convex and it can be re-written
  as
  \[ \mathcal T_{c,\epsilon}(\mu,\nu):=\inf\left\{\epsilon \mathcal H(\gamma|\bar{\gamma})\;|\; \gamma\in\Pi(\mu,\nu)\right\}, \]
  where $\mathcal H(\gamma|\bar{\gamma}):=\int_{X\times Y}(\log(\frac{\gamma}{\bar\gamma})-1)d\gamma$ is the relative entropy and $\bar{\gamma}=\exp(\frac{-c}{\epsilon})$.
  It can be proved that \eqref{eq: entropic OT} is strictly convex problem and it admits a solution of the form
   \[\gamma^\star=\exp\Big(\frac{u^\star}{\epsilon}\Big)\bar{\gamma}\exp\Big(\frac{v^\star}{\epsilon}\Big), \]
   where $u^\star$ and $v^\star$ are the Lagrange multipliers associated to the marginal constraints.
   We also highlight that by adding the entropic term we have penalized the non-negative constraint meaning that the solution $\gamma^\star$ is always positive.\\
 The regularized version of \eqref{pbVariation} can be stated in the following way
 \begin{equation}
 \label{eq:reg CN}
 \inf\{ \epsilon \mathcal H(\gamma|\bar{\gamma}) +\mathcal F(\pi_y(\gamma))+\mathcal G(\pi_x(\gamma)) \}
 \end{equation} 
where $\mathcal F$ is the functional capturing the congestion and interactions effects and $\mathcal G(\rho)=i_\mu(\rho)$ is the indicator function in the convex analysis sense and it is used to enforce the prescribed first marginal.
Before introducing the generalization of Sinkhorn algorithm proposed in \cite{Chizat} we briefly recall without proof a classical duality result
\begin{proposition}
The dual problem of \eqref{eq:reg CN}
\begin{equation}
\label{eq:duality}
\sup_{u,v}-\mathcal G^\star(-u)-\mathcal F^\star(-v)-\epsilon\int_{X\times Y}\exp\Big(\frac{u^\star}{\epsilon}\Big)\bar{\gamma}\exp\Big(\frac{v^\star}{\epsilon}\Big)dxdy.
\end{equation}
Moreover strong duality holds.
\end{proposition}
The generalized Sinkhorn algorithm is then obtained by relaxation of the maximizations on the dual problem \eqref{eq:duality}. We get the iterative method computing a sequence of potentials: given $2$ vectors $u^{(0)}$ and $v^{(0)}$, then the update at step $n$ is defined as
\begin{equation}
\label{Sinkhorn}
\begin{cases}
&u^{(n)}:=\argmax_u-\mathcal G^\star(-u)-\epsilon\int_{X\times Y}\exp\Big(\frac{u}{\epsilon}\Big)\bar{\gamma}\exp\Big(\frac{v^{(n-1)}}{\epsilon}\Big),\\
&v^{(n)}:=\argmax_v-\mathcal F^\star(-v)-\epsilon\int_{X\times Y}\exp\Big(\frac{u^{(n)}}{\epsilon}\Big)\bar{\gamma}\exp\Big(\frac{v}{\epsilon}\Big).
\end{cases}
\end{equation}
\begin{remark}
For many interesting functionals $\mathcal F$ the relaxed maximizations can be computed point-wise in space and analytically. Notice that a problem can arise in treating the interaction term, but one can easily overcome this difficulty by adopting a semi-implicit approach as the one proposed in \cite{Blanchet2017}.
\end{remark}
\paragraph{Numerical results}
$\mu$ is the uniform measure on $X=\{x_1,x_2>0\;:\; x^2_1+x^2_2<1\}$. Moreover, we take here a functional $\mathcal F$, as in the best reply scheme section, given by the sum of  a potential $V(y)=|y-\frac{\pi}{12}|^2$ and a quadratic interaction $\phi(y,z)=|y-z|^2$. 
We can then compare the numerical results obtained by the best reply scheme and 
generalized Sinkhorn, with a regularization parameter $\epsilon = 10^{-3}$, in the case of the quadratic cost function.

\begin{figure}[h!]
\TabTwo{
\includegraphics[width=.55\linewidth]{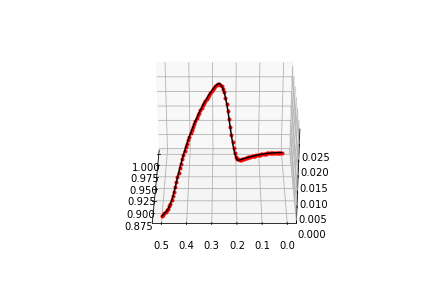}&
\includegraphics[width=.55\linewidth]{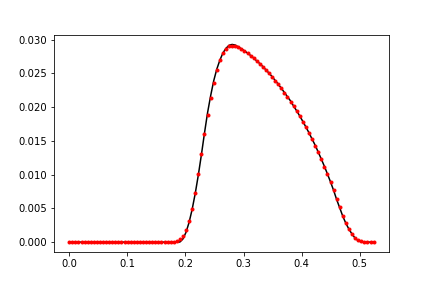}\\
}
\caption{(Left) Final density $\nu^\star$ obtained by the best reply scheme (dark solid line) and by entropic regularization (red dots) concentrated on the arc $(0,\frac{\pi}{6})$. (Right) Final density $\nu^\star$ obtained by the best reply scheme (dark solid line) and by entropic regularization (red dots).}
\label{fig:quadBRvsSinkhorn}
\end{figure}
\begin{remark}
Notice that we have not compared the iterative method for the congestion case with the entropic regularization. This is almost due to the fact that adding an other entropy term to the problem will cause an extra diffusion which could be difficult to treat even by choosing a small regularization parameter.
\end{remark}

\bibliographystyle{plain}

\bibliography{bibli}

\begin{thebibliography}{10}

\bibitem{Aumann2}
R.~Aumann.
\newblock Existence of competitive equilibria in markets with a continuum of
  traders.
\newblock {\em Econometrica}, 32:39--50, 1964.

\bibitem{Aumann}
R.~Aumann.
\newblock Markets with a continuum of traders.
\newblock {\em Econometrica}, 34:1--17, 1966.

\bibitem{benamouetalentropic}
Jean-David Benamou, Guillaume Carlier, Marco Cuturi, Luca Nenna, and Gabriel
  Peyr{\'e}.
\newblock Iterative {B}regman projections for regularized transportation
  problems.
\newblock {\em SIAM J. Sci. Comput.}, 37(2):A1111--A1138, 2015.

\bibitem{abgcptrl}
Adrien Blanchet and Guillaume Carlier.
\newblock From {N}ash to {C}ournot-{N}ash equilibria via the
  {M}onge-{K}antorovich problem.
\newblock {\em Philos. Trans. R. Soc. Lond. Ser. A Math. Phys. Eng. Sci.},
  372(2028):20130398, 11, 2014.

\bibitem{blanchet2014remarks}
Adrien Blanchet and Guillaume Carlier.
\newblock Remarks on existence and uniqueness of {C}ournot--{N}ash equilibria
  in the non-potential case.
\newblock {\em Mathematics and Financial Economics}, 8(4):417--433, 2014.

\bibitem{abgcmor}
Adrien Blanchet and Guillaume Carlier.
\newblock Optimal transport and {C}ournot-{N}ash equilibria.
\newblock {\em Math. Oper. Res.}, 41(1):125--145, 2016.

\bibitem{Blanchet2017}
Adrien Blanchet, Guillaume Carlier, and Luca Nenna.
\newblock Computation of cournot--nash equilibria by entropic regularization.
\newblock {\em Vietnam Journal of Mathematics}, Sep 2017.

\bibitem{MultiMatch}
P.-A. {Chiappori}, R.~{McCann}, and B.~{Pass}.
\newblock {Multidimensional matching}.
\newblock {\em ArXiv e-prints}, April 2016.

\bibitem{PassM2one}
Pierre-Andr{\'e} Chiappori, Robert~J. McCann, and Brendan Pass.
\newblock Multi-to one-dimensional optimal transport.
\newblock {\em Communications on Pure and Applied Mathematics}, pages n/a--n/a.

\bibitem{Chizat}
L.~Chizat, G.~Peyr\'e, B.~Schmitzer, and F.-X. Vialard.
\newblock Scaling algorithms for unbalanced transport problems.
\newblock Technical report, http://arxiv.org/abs/1607.05816, 2016.

\bibitem{CuturiSinkhorn}
M.~Cuturi.
\newblock Sinkhorn distances: Lightspeed computation of optimal transport.
\newblock In {\em Advances in Neural Information Processing Systems (NIPS) 26},
  pages 2292--2300, 2013.

\bibitem{Galichon-Entropic}
A.~Galichon and B.~Salani\'e.
\newblock Matching with trade-offs: Revealed preferences over competing
  characteristics.
\newblock Technical report, Preprint SSRN-1487307, 2009.

\bibitem{MasColell}
A.~Mas-Colell.
\newblock On a theorem of {S}chmeidler.
\newblock {\em J. Math. Econ.}, 3:201--206, 1984.

\bibitem{mccann2018optimal}
Robert~J McCann and Brendan Pass.
\newblock Optimal transportation between unequal dimensions.
\newblock {\em arXiv preprint arXiv:1805.11187}, 2018.

\bibitem{NennaPass1}
Luca Nenna and Brendan Pass.
\newblock Variational problems involving unequal dimensional optimal transport.
\newblock {\em Journal de Math{\'e}matiques Pures et Appliqu{\'e}es},
  139:83--108, 2020.

\bibitem{NennaPass3}
Luca Nenna and Brendan Pass.
\newblock Transport type metrics on the space of probability measures involving
  singular base measures.
\newblock {\em arXiv preprint arXiv:2201.00875}, 2022.

\bibitem{peyregradientflows}
Gabriel Peyr{\'e}.
\newblock Entropic approximation of {W}asserstein gradient flows.
\newblock {\em SIAM J. Imaging Sci.}, 8(4):2323--2351, 2015.

\bibitem{santambook}
Filippo Santambrogio.
\newblock {\em Optimal transport for applied mathematicians}.
\newblock Progress in Nonlinear Differential Equations and their Applications,
  87. Birkh\"auser/Springer, Cham, 2015.
\newblock Calculus of variations, PDEs, and modeling.

\bibitem{Villani-OptimalTransport-09}
C.~Villani.
\newblock {\em Optimal Transport: Old and New}, volume 338 of {\em Grundlehren
  der mathematischen Wissenschaften}.
\newblock Springer, 2009.

\end{thebibliography}

\end{document}